\documentclass[a4paper,12pt]{amsart}
\usepackage[all]{xy}
\usepackage{wrapfig}
\usepackage[dvips]{graphicx}
\usepackage{amsmath,amsthm,amssymb}

\theoremstyle{definition}
\newtheorem{thm}{Theorem}[section]
\theoremstyle{definition}

\theoremstyle{definition}
\newtheorem{lem}[thm]{Lemma}
\theoremstyle{definition}
\newtheorem{df}[thm]{Definition}
\theoremstyle{definition}
\newtheorem{prop}[thm]{Proposition}
\theoremstyle{definition}
\newtheorem{ex}[thm]{Example}
\theoremstyle{definition}
\newtheorem{rem}[thm]{Remark}
\theoremstyle{definition}

\theoremstyle{definition}

\usepackage{amscd}

\makeindex
\title{Cuntz's $ax+b$-semigroup $C^*$-algebra over $\mathbb{N}$ and product system $C^*$-algebras}
\author{Shinji Yamashita}
\address{Graduate School of Mathmatics, Kyushu University,
Hakozaki, Fukuoka, 812-8581, JAPAN}
\email{s-yamashita@math.kyushu-u.ac.jp}
\date{\today}

\begin{document}
\maketitle

\begin{abstract}
We investigate $C^*$-algebras associated with row-finite topological higher-rank graphs with no source, which are based on product system $C^*$-algebras. We prove the Cuntz-Krieger uniqueness theorem, and provide the condition of simplicity and purely infiniteness of our algebras. We give examples of non-discrete topological higher-rank graphs whose $C^*$-algebras contain Cuntz's $ax+b$-semigroup $C^*$-algebra over $\mathbb{N}$.
\end{abstract}


\section{Introduction}
\label{sec:Introduction}

Recently, Cuntz introduced a $C^*$-algebra $\mathcal{Q}_\mathbb{N}$ which is naturally associated with the $ax+b$-semigroup over $\mathbb{N}$ (\cite{Cuntz}) and this algebra seems some kind of quantization of natural numbers. He showed that it is simple and purely infinite and generated by Bost-Connes $C^*$-algebra $\mathcal{C}_\mathbb{Q}$ (\cite{BostConnes}) and one unitary. Moreover, Cuntz and Li extended this construction to some commutative rings without zero divisors (\cite{CuntzLi}) which is motivated the extension of the Bost-Connes algebra to arbitrary number fields, and recently, Li has done for arbitrary rings (\cite{Li}). Moreover, the $C^*$-algebra $\mathcal{Q}_\mathbb{N}$ has an arithmetic flavor from the construction and we think this $C^*$-algebra will become a brigde between operator algebras and number theory like the Bost-Connes $C^*$-algebra. 

Including such algebras, some important properties of $C^*$-algebras are determined by the dynamical system on the underlying space. In this fashion, Katsura introduced a new class of $C^*$-algebras which is called a topological graph $C^*$-algebra (\cite{Katsura1}). This class is not only a generalization of both graph algebras and homeomorphism $C^*$-algebras, but also many interesting $C^*$-algebras, for example, all Kirchberg $C^*$-algebras are included in this class (\cite{Katsura4}), and he systematically studied these algebras in a series of his paper.

On the other hand, Kumjian and Pask (\cite{KumjianPask}) introduced a higher-dimensional graph $C^*$-algebra which is called a higher-rank graph $C^*$-algebra and Raeburn, Sims and Yeend extended it which is now called a finitely aligned higher-rank $C^*$-algebra (\cite{RaeburnSimsYeend1}, \cite{RaeburnSimsYeend2}). One of the most interesting and attractive facts which does not observe in graph algebras is to appear A$\mathbb{T}$-algebras like Bunce-Deddens algebras and irrational rotation algebras as a higher-rank graph $C^*$-algebra (\cite{PaskRaeburnRordamSims}). We feel that higher-rank $C^*$-algebras have fruitful structures and are not merely higher dimensional version of graph $C^*$-algebras.

Under these backgrounds, in \cite{Yeend1}, \cite{Yeend2}, Yeend introduced a notion of a topological higher-rank graph which is unify a higher-rank graph and a topological graph $C^*$-algebra. His approach is based on the theory of the groupoid $C^*$-algebras (\cite{Renault}). On the other hand, as one of higher-dimensional models of a Cuntz-Pimsner algebra, Fowler considered product system $C^*$-algebras associated with quasi-lattice ordered group $(G,P)$ in \cite{Fowler} and Raeburn and Sims researched relations between product system $C^*$-algebras and higher-rank graph $C^*$-algebras in \cite{RaeburnSims}.

 The purpose of this paper is to study a 'row-finite with no source' topological higher-rank graph based on product system $C^*$-algebras.
 In particular, we establish a method of an analysis of a (non-discrete) topological higher-rank graph and construct a new class of $C^*$-algebras which can be considered as some extensions of Cuntz's $ax+b$-semigroup $C^*$-algebras. 
 It should be mentioned that Yeend commented in the introduction of \cite{Yeend1} that there is a problem in product system $C^*$-algebras, in particular the definition of Cuntz-Pimsner covariance, and in later, Sims and Yeend (\cite{SimsYeend}) suggest a new Cuntz-Pimsner covariance of a product system $C^*$-algebra which includes $C^*$-correspondences (\cite{Katsura2}), finitely aligned higher-rank $C^*$-algebras (\cite{RaeburnSimsYeend2}) and coincide with Fowler's one under the 'row-finite with no source' case. 
 However, there is no research to establish a theory of the $C^*$-dynamical system
 based on the product system $C^*$-algebras and topological higher-rank graphs. Furthermore we need slightly different techniques compared with (discrete) higher-rank graphs.
 
 In this paper, we actively consider the infinite dimension case of topological higher-rank graph though the dimension of higher-rank graphs is supposed to be finite in usual. In terms of product system $C^*$-algebras, it is equivalent to consider the ordered lattice group $(\oplus_{i=1}^\infty \mathbb{Z} , \oplus_{i=1}^\infty \mathbb{N} )$. One of this reason is that $\mathcal{Q}_\mathbb{N}$ (and also $\mathcal{C}_\mathbb{Q}$) have a structure based on the prime numbers (more precisely, the prime factorization of natural numbers) which is corresponding to the basis of $(\oplus_{i=1}^\infty \mathbb{Z} , \oplus_{i=1}^\infty \mathbb{N} )$. We can expect that this case leads more fruitful structures for topological higher-rank graph $C^*$-algebras.
 
 This paper is organised as follows. In Section 2, we define a topological higher-rank graph followed by Yeend and recall several definitions and fundamental facts.
 In Section 3, we prove the Cuntz-Krieger Uniqueness theorem under the assumption 'row-finite with no source' (Theorem \ref{thm:CK-unique}). This is our main theorem and a generalization of Theorem 4.6 of \cite{KumjianPask}. A point of the proof of this theorem is to overcome the difficulty of calculations of some multiplications which does not occur in higher-rank graph $C^*$-algebras.
 In Section 4, we give the condition for topological higher-rank $C^*$-algebras to be simple and purely infinite.
 In Section 5, we construct a higher rank version of a topological graph investigated by Deaconu (\cite{Deaconu1}) and Katsura (\cite{Katsura4}). The associated $C^*$-algebras include Cuntz's $ax+b$-semigroup $C^*$-algebra $\mathcal{Q}_\mathbb{N}$ over $\mathbb{N}$. The author think that this way of the extension of the $ax+b$-semigroup $C^*$-algebra is different to the arithmetic one noted as above, but we expect some relationships with number theory.
 
After completed this work, the author has found that Renault, Sims and Yeend proved the uniqueness theorem for a $C^*$-algebra associated with a compactly aligned topological higher-rank graph which is more generalized setting from the author's one (\cite{RenaultSimsYeend}). However they proved this theorem via a groupoid theory (\cite{Renault}), on the other hand, the author's method is used a theory of a product system $C^*$-algebra.  
 
\section{Preliminaries}
\label{sec:Preliminaries}

In this section, we recall a topological $k$-graph defined by Yeend and a product system $C^*$-algebra defined by Fowler, in order to construct an associated $C^*$-algebra for a topological $k$-graph.
 
 We denote the set of natural numbers by $\mathbb{N}=\{ 0,1,2,\cdots \}$ and the integers by $\mathbb{Z}$.   We denote by $\mathbb{T}$ the group consisting of complex numbers whose absolute values are 1. Given a (semi)group with identity $P$ and $k=1,2,\cdots, \infty$, we denote $P^k=\otimes_{i=1}^k P$ by the direct sum of $P$ which has a natural (semi)group structure.
 We consider $\mathbb{N}^k$ as an additive semigroup with identity 0.
 For $1 \le k \le \infty$, $e_1,e_2,\cdots$ are standard generators of $\mathbb{N}^k$.
 We write the $i$-th coordinate of $m \in \mathbb{N}^k$ by $m_{(i)}$.
 For $m,n \in \mathbb{N}^k$, we say $m \le n$ if $m_{(i)} \le n_{(i)}$ for any $1 \le i \le k$.
 For a locally compact (Hausdorff) space $\Omega$, we denote by $C(\Omega )$ the linear space of all continuous functions on $\Omega$. We define $C_c(\Omega )$, $C_0(\Omega )$, $C_b(\Omega )$ by those of compactly supported functions, functions vanishing at infinity, and bounded functions, respectively. 
 
Next we shall define a topological $k$-graph which was defined by Yeend (\cite{Yeend1}, \cite{Yeend2}).
\begin{df}
Given $k=1,2,\cdots ,\infty$, a \textit{topological $k$-graph} is a pair $(\Lambda , d)$ consisting of a small category $\Lambda = ( \mathrm{Obj}(\Lambda ) , \mathrm{Mor}(\Lambda ) , r,s)$ and a functor $d:\Lambda \longrightarrow \mathbb{N}^k$, called the \textit{degree map}, which satisfy the following:
\begin{enumerate}
 \item{$\mathrm{Obj}(\Lambda )$ and $\mathrm{Mor}(\Lambda )$ are locally compact spaces;}
 \item{$r,s: \mathrm{Mor}(\Lambda ) \longrightarrow \mathrm{Obj}(\Lambda )$ are continuous and $s$ is a local homeomorphism;}
 \item{Composition $\circ : \Lambda \times_c \Lambda \longrightarrow \Lambda$ is continuous and open, where $\Lambda \times_c \Lambda=\{ (\lambda , \mu ) | s(\lambda )=r(\mu ) \}$ has the relative topology inherited from the product topology on $\Lambda \times \Lambda$;}
 \item{$d$ is continuous, where $\mathbb{N}^k$ has the discrete topology;}
 \item{For all $\lambda \in \Lambda$ and $m,n \in \mathbb{N}^k$ such that $d(\lambda )=m+n$, there exists unique $(\mu, \nu )\in \Lambda \times_c \Lambda$ such that $\lambda = \mu \nu , d(\mu )=m$ and $d(\nu )=n$}.
\end{enumerate} 
\end{df}

In \cite{Yeend1} and \cite{Yeend2}, Yeend supposed the second-countablility of the locally compact spaces $\mathrm{Obj}(\Lambda )$ and $\mathrm{Mor}(\Lambda )$. However we will not use this assumption in here.
 We refer to the morphisms of $\Lambda$ as \textit{paths} and to the objects of $\Lambda$ as \textit{vertices}.
 The maps $r$ and $s$ are called the \textit{range} and \textit{source maps}, respectively.
 Given subsets $U,V \subset \Lambda$, we define $U \times_c V=\{ (\lambda , \mu ) \in U \times V | s(\lambda )=r(\mu ) \}$

 We give an example of topological $k$-graph.
 Set $\mathrm{Obj}(\Omega_k)=\mathbb{N}^k$ and $\mathrm{Mor}(\Omega_k)=\{ (m,n) \in \mathbb{N}^k \times \mathbb{N}^k | m \le n \}$ with discrete topologies.
 The range and source maps are given by $r(m,n)=m$, $s(m,n)=n$.
 Define $d: \Omega_k \longrightarrow \mathbb{N}^k$ by $d(m,n)=n-m$.
 Then $\Omega_k$ is a topological $k$-graph. 

 For $m \in \mathbb{N}^k$, define $\Lambda^m$ to be the set $d^{-1}(\{ m \})$ of paths of degree $m$.
 Then $\Lambda^m \times_c \Lambda^n$ and $\Lambda^{m+n}$ are homeomorphic by the restriction of the composition map.
 Let us put $r_m , s_m : \Lambda^m \longrightarrow \Lambda^0$ by the restriction maps of $r,s$ on $\Lambda^m$.
 Then $(\Lambda^0 , \Lambda^m, s_m , r_m )$ is a topological graph defined in Katsura \cite{Katsura1}.
 Set $\Lambda^m(v)=r_m^{-1}(v)$ for $v \in \Lambda^0$.
 For $\lambda \in \Lambda^l$ and $0 \le m \le n \le l$, there exists uniquely $\lambda_1 \in \Lambda^m$, $\lambda_2 \in \Lambda^{n-m}$ and $\lambda_3 \in \Lambda^{l-n}$ such that $\lambda = \lambda_1 \lambda_2 \lambda_3$ by the factorization property.
 We shall denote $\lambda_2$ as $\lambda (m,n)$, and especially, we shall write $\lambda (m)=\lambda (m,m)$.
 Then we can define the continuous map $\mathrm{Seg}_{(m,n)}^l : \Lambda^l \longrightarrow \Lambda^{n-m}$ by $\mathrm{Seg}_{(m,n)}^l(\lambda)=\lambda(m,n)$.

 Let $(\Lambda_1,d_1)$ and $(\Lambda_2,d_2)$ be topological $k$-graphs.
 A \textit{$k$-graph morphism} between $\Lambda_1$ and $\Lambda_2$ is a continuous functor $\alpha : \Lambda_1 \longrightarrow \Lambda_2$ satisfying $d_2(\alpha (\lambda ))=d_1(\lambda )$ for all $\lambda \in \Lambda_1$.

Next we shall define a property of $\Lambda$ which is correspond to Definition of 1.4 of \cite{KumjianPask} for a higher-rank graph. 
\begin{df}
The topological $k$-graph $\Lambda$ is \textit{row-finite for degree $m$} $(m \in \mathbb{N}^k)$ if for each $v \in \Lambda^0$, there exists a neighborhood $V \subset \Lambda^m$ of $v$ such that $r_m^{-1}(V)$ is a compact set in $\Lambda^m$. $\Lambda$ \textit{has no source for degree $m$} for any $v \in \Lambda^0$ and any neighborhood $V \subset \Lambda^m$ of $v$, $r_m^{-1}(V)\neq \emptyset$. We say $\Lambda$ is \textit{row-finite} (resp. \textit{has no source}) if for every $m \in \mathbb{N}^k$, $\Lambda$ is row-finite (resp. has no source) for degree $m$.    
\end{df}
 
If $\Lambda^{e_i}$ is row-finite (resp. has no source) for every standard generators $e_i \in \mathbb{N}^k$, then $\Lambda$ is row-finite (resp. has no source) by the following proposition.

\begin{prop}
\label{prop:test}
\normalfont\slshape
 For $m,n \in \mathbb{N}^k$,
$\Lambda$ are row-finite (resp. has no source) for degree $m,n$, then $\Lambda$ is row-finite (resp. has no source) for degree $m+n$. 
\end{prop}

\begin{proof}
First, we suppose that $\Lambda^m, \Lambda^n$ are row-finite. For any $v \in \Lambda^m$, $r_m^{-1}(v)$ is compact, hence $s_m(r_m^{-1}(v))$ is also compact. Since $\Lambda^n$ is row-finite, we can take a set of neighborhoods $\{ W_1,\cdots ,W_l \}$ such that
$s_m(r_m^{-1}(v))$ is covered by $\{W_i \}_{i=1}^l$ and $r_n^{-1}(W_i)$ is compact.
This implies $r_m^{-1}(v)$ is covered by $\{s_m^{-1}(W_i) \}_{i=1}^l$. By Lemma 1.21 of Katsura's article, there exists a neighborhood $V$ of $v$ such that 
\begin{equation*}
r_m^{-1}(V) \subset  \bigcup_{i=1}^l s_m^{-1}(W_i).
\end{equation*}
and $r_m^{-1}(V)$ is compact since $\Lambda$ is row-finite for degree $n$. Then
\begin{equation*}
r_{m+n}^{-1}(V)=r_m^{-1}(U) \times_c \bigcup_{i=1}^l r_n^{-1}(W_i)
\end{equation*}
and this implies $\Lambda$ is row-finite for degree $m+n$.

On the other hand, since $\Lambda^m$ has no source, for any $v \in \Lambda^0$, there is a neighborhood $V$ of $v$ such that $r_m^{-1}(V) \neq \emptyset$. Hence
\begin{equation*}
r_{m+n}^{-1}(V)=r_m^{-1}(V) \times_c r_n^{-1}(s_m(r_m^{-1}(V))) \neq \emptyset ,
\end{equation*}
hence $\Lambda^{m+n}$ has no source.
\end{proof}

Let $\Lambda$ be a row-finite topological $k$-graph with no source. We define the \textit{infinite path space} of $\Lambda$ by
\begin{center}
$\Lambda^\infty$ = $\{ \alpha : \Omega_k \longrightarrow \Lambda$ $|$ $\alpha$ is a $k$-graph morphism $\}$. 
\end{center}
We extend the range map by setting $r(\alpha )=\alpha (0)$.
 Set $\Lambda^\infty (v) = \{ \alpha \in \Lambda^\infty | v= \alpha (0 ) \}$ for $v \in \Lambda^0$.
 Remark that $\Lambda^\infty (v) \neq \emptyset$ since we assume $\Lambda$ has no source.
 For each $p \in \mathbb{N}^k$, define $\tau^p : \Lambda^\infty \longrightarrow \Lambda^\infty$ by $\tau^p(\alpha )(m,n)=\alpha (m+p, n+p)$ for $\alpha \in \Lambda^\infty$ and $(m,n) \in \Omega_k$.

 To define $C^*$-algebras from topological higher-rank graphs, we use the product system $C^*$-algebras considered by Fowler \cite{Fowler}.
 First we shall recall a Hilbert $A$-bimodule.
 Let $A$ be a $C^*$-algebra and $Y$ be a right Hilbert $A$-module with an $A$-valued inner product $\langle \cdot , \cdot \rangle$.
 We denote by $L(Y)$ the $C^*$-algebra of the adjointable operators on $Y$.
 Given $x,y \in Y$, the rank-one operator $\theta_{x,y} \in L(Y)$ is defined by $\theta_{x,y}(z)=x \langle y,z \rangle$ for $z \in Y$.
 The closure of the linear span of rank-one operators is denoted by $K(Y)$.
 We say that $Y$ is a Hilbert $A$-bimodule if $Y$ is a right Hilbert $A$-module with a homomorphism $\phi : A \longrightarrow L(Y)$.
 Then we can define a left $A$-action on $Y$ by $a \cdot x=\phi (a)x$.
 
 Let $P$ be a discrete multiplicative semigroup with identity $e$, and let $A$ be a $C^*$-algebra.
 A \textit{product system} over $P$ of Hilbert $A$-bimodules is a semigroup $X= \bigsqcup_{p \in P} X_p$ such that (1) for each $p \in P$, $X_p \subset X$ is a Hilbert $A$-bimodule with a left action $\phi_p$ of $A$ ;(2) the identity fibre $X_e$ is equal to the bimodule ${}_AA_A$; (3) for $p,q \in P \setminus \{ e \}$, there is an isomorphism $M_{p,q} : X_p \otimes_A X_q \longrightarrow X_{pq}$ (in this paper, we also denote $x \otimes y \in X_{pq}$ instead of $M_{p,q}(x \otimes y)$) ; (4) multiplication in $X$ by elements of $X_e=A$ implements the action of $A$ on each $X_p$ ; that is $M_{e,p}(a \otimes x)=a \cdot x$ and $M_{p,e}(x \otimes a)=x \cdot a$ for all $p \in P , x \in X_p$ and $a \in X_e=A$.
In this paper, we always suppose $P$ is subsemigroup of a discrete group $G$ such that $P \cap P^{-1}=\{ e\}$ and with respect to the partial order $p \le q$ $\iff$ $p^{-1}q \in P$, any two element $p,q \in G$ which have a common upper bound in $P$ have a least upper bound $p \vee q \in P$.
 We write $p \vee q = \infty$ to indicate that $p,q \in G$ have no common upper bound in $P$.
 We say the pair $(G,P)$ with such properties is a \textit{quasi-lattice ordered group}.
 The main example of interest to us is $(G,P)=(\mathbb{Z}^k, \mathbb{N}^k)$, which is actually lattice-ordered: each $m,n \in \mathbb{N}^k$ has a least upper bound $m \vee n$ with $i$-th coordinate $(m \vee n)_{(i)}= \max \{ m_{(i)} , n_{(i)} \}$.

 Let $(G,P)$ be a quasi-lattice ordered group and let $X$ be a product system over $P$.
 Throughout this paper, we suppose the left action of $A$ on each fiber $X_p$ is an injective into $K(X_p)$.  
 Let $B$ be a $C^*$-algebra, and let $T$ be a map from $X$ to $B$. 
 For $p \in P$, let $T_p=T|_{X_p}$. We say that $T$ is a \textit{representation} of $X$ if (i) $T_e$ is *-homomorphism and each $T_p$ is linear (ii) $T_p(x)T_q(y)=T_{pq}(xy)$ for $p,q \in P$ and $x \in X_p$, $y \in X_q$ (iii) $T_p(x)^*T_p(y)=T_e(\langle x,y \rangle )$ for $x,y \in X_p$. 
 We call $T$ is \textit{injective} if $T_e$ is injective.
 Let us define $C^*(T)=\overline{\mathrm{span}} \{ T_p(x)T_q(y)^* | p,q \in P , x \in X_p ,y \in X_q \}$ and $\mathcal{F}_T=\overline{\mathrm{span}} \{ T_p(x)T_p(y)^* | p \in P , x ,y \in X_p \}$ which is called the core of $C^*(T)$.
 Define $\psi_p^T: K(X_p) \longrightarrow C^*(T)$ by $\psi_p^T(\theta_{x,y} )=T_p(x)T_p(y)^*$, then $\psi_p^T$ is well-defined by Lemma 2.2 of \cite{KajiwaraPinzariWatatani}.
 We say a representation $T$ is \textit{Cuntz-Pimsner covariance} if for each $p \in P$, $\psi_p^T \circ \phi_p(a) = T_0(a)$ for any $a \in A$.
 We say that $X$ is \textit{compactly aligned} if for all $p,q \in P$ such that $p \vee q < \infty$, and $S_1 \in K(X_p)$ and $S_2 \in K(X_q)$, we have $(S_1 \otimes 1_{p^{-1}(p \vee q )})(S_2 \otimes 1_{q^{-1}(p \vee q)}) \in K(X_{p \vee q})$.
 If the image of the left action of $A$ on each fiber $X_p$ is in $K(X_p)$ which is assumed in this paper, then $X$ is compactly aligned by Proposition 5.8 of \cite{Fowler}.
 In this case, $T$ satisfies the following relation which is called \textit{Nica covariance} by Proposition 5.4 of \cite{Fowler}; for $S_1 \in K(X_p),\ S_2 \in K(X_q)$,
\[ 
\psi^T_p(S_1)\psi^T_q(S_2)=
\left\{
  \begin{array}{ll}
  \psi_{p \vee q}^T((S_1 \otimes 1_{p^{-1}(p \vee q )})(S_2 \otimes 1_{q^{-1}(p \vee q)})) & \mbox{if $p \vee q < \infty$} \\
  0 & \mbox{if $p \vee q =\infty$}
  \end{array}
\right.
\]

 There is a $C^*$-algebra $\mathcal{O}_X$ and a Cuntz-Pimsner covariant representation $t:X \longrightarrow \mathcal{O}_X$ which is universal in the following sense: $\mathcal{O}_X=C^*(t)$ and for any Cuntz-Pimsner covariance $T: X \longrightarrow B$, there is a unique homomorphism $T_* : \mathcal{O}_X \longrightarrow B$ such that $T_* \circ t = T$. Remark that $t$ is injective (see \cite{Fowler}).

 For $P=\mathbb{N}^k$, the universality allows us to define a strongly continuous gauge action $\gamma : \widehat{\mathbb{Z}^k}=\prod_{i=1}^k \mathbb{T} \longrightarrow \mathrm{Aut}(\mathcal{O}_X)$ such that $\gamma_z(t_n(x))=z^n(t_n(x))$ for $z \in \prod_{i=1}^k \mathbb{T}$, $x \in X_n$ and the fixed point algebra $\mathcal{O}_X^\gamma$ is $\mathcal{F}_t=\overline{\mathrm{span}}\{T_n(x)T_n(y)^* | n \in \mathbb{N}^k, x,y \in X_n \}$. Let
\begin{equation*}
\Psi(x)= \int_{z \in \prod_{i=1}^k \mathbb{T}} \gamma_z (x)dz \quad (x \in \mathcal{O}_X)
\end{equation*}
be a faithful conditional expectation onto $\mathcal{O}_X^\gamma$.

 We shall consider a product system $C^*$-algebra associated with a row-finite topological $k$-graph $\Lambda$ with no source. For the quasi-lattice ordered group $(G,P)=(\mathbb{Z}^k, \mathbb{N}^k)$, define a $C^*$-algebra $A=C_0(\Lambda^0)$ and 
\begin{equation*}
X_n=C_s(\Lambda^n)=\{ \xi \in C(\Lambda^n) | \langle \xi, \xi \rangle \in C_0(\Lambda^0) \}
\end{equation*}
where the inner product $\langle \cdot , \cdot \rangle$ is defined by
\begin{equation*}
\langle \xi ,\eta \rangle (v)= \sum_{s_n(\lambda )=v } \overline{\xi (\lambda )} \eta (\lambda )
\end{equation*} 
for $\xi,\eta \in X_n$ and $v \in \Lambda^0$, and the left and right actions are defined by
\begin{equation*}
(f\cdot \xi \cdot g)(\lambda )=f(r_n(\lambda )) \xi (\lambda ) g(s_n(\lambda ))
\end{equation*}
for $f,g \in A$, $\xi \in X_n$, $n \in \mathbb{N}^k$. Define $M_{m,n}:X_m \otimes_A X_n \longrightarrow X_{m+n}$ by $M_{m,n}(\xi \otimes \eta)(\lambda )=\xi(\lambda (0,m)) \eta(\lambda (m,m+n))$ for $m,n \in \mathbb{N}^k$. Then we can construct a product system $X$ from these data.

\begin{df}
Let $\Lambda$ be a row-finite topological $k$-graph with no source and $X$ is a product system defined as above. We define the $C^*$-algebra $\mathcal{O}(\Lambda )$ to be the universal $C^*$-algebra $\mathcal{O}_X$ built from the product system $X$.
\end{df}

 For a product system $X$ associated with a row-finite topological $k$-graph with no source and a Cuntz-Pimsner covariance $T:X \longrightarrow C^*(T)$, let us define a sub *-algebra $C^*(T)^\mathrm{cpt}$ of $C^*(T)$ by
\begin{equation*}
C^*(T)^\mathrm{cpt}=\mathrm{span} \{ T_n(\xi ) T_m(\eta )^* | n,m \in \mathbb{N}^k,\ \xi \in C_c(\Lambda^n ) ,\ \eta \in C_c(\Lambda^m) \} .
\end{equation*}
 This sub *-algebra is dense in $C^*(T)$ by Lemma 1.6 of \cite{Katsura1}. We will use the dense algebra $C^*(T)^\mathrm{cpt}$ instead of the algebraic part of $C^*(T)$.
 

\section{Cuntz-Krieger Uniqueness Theorem}
\label{sec:CK-Unique}

In this section, we shall prove so-called Cuntz-Kriger Uniqueness Theorem for certain row-finite topological $k$-graph with no source (Theorem \ref{thm:CK-unique}).
 The phenomena which does not occur in higher-rank graph $C^*$-algebras or topological (1-)graph $C^*$-algebras is to be difficult to compute the form $T_n(\xi )^*T_m(\eta )$ using the terms of the underlying topological $k$-graph if $n$ and $m$ are unordered (see Lemma 3.1 of \cite{KumjianPask} and Lemma 2.4 of \cite{Katsura1}).
 To avoid this problem, we prove this theorem by the combination of the idea of \cite{KajiwaraPinzariWatatani} for Cuntz-Pimsner algebras $\mathcal{O}_Y$ for Hilbert $A$-bimodule $Y$ and the notion of the Nica covariance.
 However, Kajiwara-Pinzari-Watatani (\cite{KajiwaraPinzariWatatani}) considered in the case that $Y$ has a finite basis which is stronger than our assumption, we need to modify a little bit.

First, we introduce some terms for working smoothly.

\begin{df}
Given finitely many functions $\xi_1,\cdots ,\xi_L , \eta_1 ,\cdots \eta_L$ of $C_c(\Lambda^m)$, we say $\{ ( \xi_i,\eta_i ) \}_{i=1}^L$ is an \textit{orthogonal pair} for degree $m$ if for any $i=1,\cdots ,L$, $\xi_i(\lambda ) \overline{\eta_i (\lambda ' )}=0$ for $s(\lambda )=s(\lambda ')$ and $\lambda \neq \lambda '$.
 For a set $\Omega \subset \Lambda^m$ and $u_1,\cdots ,u_L \in C_c(\Lambda^m)$, $\{u_i\}_{i=1}^L$ is a \textit{partition of unity} for $\Omega$ if $u_i$ satisfies $0 \le u_i \le 1$, $\sum_{i=1}^L u_i^2(\lambda )=1$ for $\lambda \in \Omega$, and $s_m$ is injective on the support $\mathrm{supp}(u_i)$ of $u_i$. In partcular, $\{ (u_i,u_i) \}_{i=1}^L$ is an orthogonal pair.
\end{df}

\begin{lem}
\normalfont\slshape
\label{lem:unity}
Let $\Lambda$ be a row-finite topological $k$-graph with no source.
\begin{enumerate}
 \item{For $\xi \in C_c(\Lambda^m)$, there exists a partition of unity $\{u_i\}_{i=1}^L$ for the compact support $\mathrm{supp}(\xi )$ of $\xi$.}
 \item{If $\xi \in X_m$ and $\{ u_i\}_{i=1}^L$ is a partition of unity for the support $\mathrm{supp}(\xi )$ of $\xi $, then $\xi=\sum_{i=1}^L u_i \langle u_i , \xi \rangle $.}
\end{enumerate}
\end{lem}

\begin{proof}
(1) Take $\xi \in C_c(\Lambda^m)$. Since $s$ is a local homeomorphism, for each $\lambda \in \Lambda^m$ there exists an relative compact open neighborhood $U_\lambda$ of $\lambda$ such that the restriction of $s_m$ to $U_\lambda$ is injective. Since $\mathrm{supp}(\xi )$ is compact, we can find $\lambda_1 ,\cdots ,\lambda_L$ such that $\mathrm{supp}(\xi ) \subset \cup_{i=1}^L U_{\lambda_i}$. Take $v_1,\cdots , v_L$ satisfying $0 \le v_i \le 1$, $\mathrm{supp}(v_i) \subset U_{\lambda_i}$ for each $1 \le i \le L$, and $\sum_{i=1}^L v_i(\lambda )=1$ for all $\lambda \in \mathrm{supp}(\xi )$. Then $u_i=v_i^{1/2}$ is a partition of unity for $\mathrm{supp}(\xi )$.\\
(2) Take $\xi \in X_m$ and a partition of unity $\{ u_i\}_{i=1}^L$ for $\mathrm{supp}(\xi )$. For $\lambda \in \Lambda^m$,
\begin{equation*}
\Bigl( \sum_{i=1}^L u_i \langle u_i , \xi \rangle \Bigr) (\lambda )
=\sum_{i=1}^L u_i(\lambda ) \Bigl( \sum_{s(\mu )=s(\lambda )}u_i(\mu ) \xi (\mu )  \Bigr)
=\sum_{i=1}^L u_i(\lambda )^2 \xi (\lambda ) = \xi (\lambda ). 
\end{equation*}
\end{proof}

\begin{lem}
\label{lem:construct unity}
\normalfont\slshape
 Let $\Lambda$ be a row-finite topological $k$-graph with no source.
 For $i=1,\cdots ,L$, let $K_i$ be a compact set of $\Lambda^{n_i}$ and put $n=\vee_{i=1}^L n_i$, $m_i=n-n_i$.
 Then, for each $i$, there are partitions of unity $\{ v_{i,q} \}_q$ for $K_i$ and $\{u_{i,p}\}_p$ for $\mathrm{Seg}_{(n_i,n)}^n(K_i \times_c \Lambda^{m_i})$. Furthermore, for each $i$, $\{ v_{i,q} \otimes u_{i,p} \}_{p,q}$ is a partition of unity for $\cup_{i=1}^L (K_i \times_c \Lambda^{m_i})$.
\end{lem}

\begin{proof}
Since we assume $\Lambda$ is row-finite, $K_i \times_c \Lambda^{m_i}$ is compact subset of $\Lambda^n$. By Lemma \ref{lem:unity}, there is a partition unity $\{w_p\}_p$ for $\cup_{i=1}^L (K_i \times_c \Lambda^{m_i})$.
 For each $w_p$, define $u_{i,p} \in C_c(\mathrm{Seg}_{(n_i,n)}^n(\mathrm{supp}(w_p))$ by the following:
 Given $\mu \in \Lambda^{m_i}$, if there exists an element $\lambda \in \Lambda^{n_i}$ such that $s_{n_i}(\lambda )= r_{m_i} (\mu )$ and $\lambda \mu \in \mathrm{supp}(w_p)$ (in this case, $\lambda$ is uniquely determined), then we define $u_{i,p}(\mu )=w_p(\lambda \mu )$ and otherwise $u_{i,p}(\mu )=0$.
 From the construction of $u_{i,p}$, $\mathrm{supp}(u_{i,p})=\mathrm{Seg}_{(n_i,n)}^n(\mathrm{supp}(w_p))$ and $\{ u_{i,p} \}_p$ is a partition of unity for $\mathrm{Seg}_{(n_i,n)}^n(K_i \times_c \Lambda^{m_i})$.
 On the other hand, there is a partition of unity $\{ v_{i,q} \}_q$ for the compact set $K_i$. Then $\{ v_{i,q} \otimes u_{i,p} \}_{p,q}$ is a partition of unity for $\cup_{i=1}^L K_i \times_c \Lambda^{m_i}$ since for each $p$, 
\begin{equation*}
\mathrm{supp}(w_p)
=\mathrm{Seg}_{(0,n_i)}^n (\mathrm{supp}(w_p)) \times_c \mathrm{Seg}_{(n_i,n)}^n (\mathrm{supp}(w_p) )
\subset \big( \bigcup_q \mathrm{supp}(v_{i,q}) \big) \times_c \mathrm{supp}(u_{i,p}).
\end{equation*}
\end{proof}

Given finitely many functions $\{ u_i\}_{i=1}^L$ in $X_m$, let us define a positive linear map $\sigma_{\{u_i\}_{i=1}^L}^{T,m} : C^*(T) \longrightarrow  C^*(T)$ by 
\begin{equation*}
\sigma_{\{u_i\}_{i=1}^L}^{T,m}(x)=\sum_{i=1}^L T_m(u_i)xT_m(u_i)^* ,\quad x \in C^*(T)
\end{equation*}

\begin{lem}
\label{lem:turn1}
\normalfont\slshape
Let $\{u_i\}_{i=1}^L$ be a partition of unity for the support $\mathrm{supp}(\xi)$ of $\xi \in X_n$. Then
$
\sigma_{\{u_i\}_i}^{T,n}(x)T_n(\xi )=T_n(\xi ) x
$
holds for every element $x$ of the relative commutant algebra $T_0(A)' \cap C^*(T)$.
\end{lem}

\begin{proof}
For $x \in T_0(A)' \cap C^*(T)$,
\begin{eqnarray*}
\sum_{i=1}^L T_n(u_i)xT_n(u_i)^*T_n(\xi )
&=&\sum_{i=1}^L T_n(u_i)xT_0(\langle u_i, \xi \rangle)
=\sum_{i=1}^L T_n(u_i)T_0(\langle u_i, \xi \rangle)x \\
&=&T_n(\sum_{i=1}^L u_i \langle u_i, \xi \rangle )x.
\end{eqnarray*}
Since $\{u_i\}_{i=1}^L$ is a partition of unity for $\mathrm{supp}(\xi)$, we obtain $\sum_{i=1}^L u_i \langle u_i, \xi \rangle =\xi$ by Lemma \ref{lem:unity}. 
\end{proof}

\begin{lem}
\label{lem:turn2}
\normalfont\slshape
For $\xi_i \in C_c(\Lambda^{m_i})$ $(i=1,\cdots ,L)$, let $K_i$ be a compact support $\mathrm{supp}(\xi_i)$.
Take a partition of unity $\{ v_{i,q} \}_q$ for $K_i$ and $u_{i,p} \in C_c(\Lambda^{n_i})$ by Lemma \ref{lem:construct unity} where $n=\vee_{i=1}^L n_i$ and $m_i=n-n_i$. Then 
\begin{equation*}
\sigma_{\{v_{i,q} \otimes u_{i,p} \}_{p,q}}^{T,n}(x)T_{n_i}(\xi_i)=T_{n_i}(\xi_i)\sigma_{\{u_{i,p}\}}^{T,m_i}(x)
\end{equation*}
holds for $x \in T_0(A)' \cap C^*(T)$.
\end{lem}

\begin{proof}
First, we consider the left action $\phi_{m_i}(\langle v_{i,q},\xi_i \rangle )$ on $X_{m_i}$. For $\zeta \in X_{m_i}$ and $\lambda \in \Lambda^{m_i}$, 
\begin{equation*}
(\phi_{m_i}(\langle v_{i,q},\xi_i \rangle )\zeta )(\lambda )
=\langle v_{i,q},\xi_i \rangle (r(\lambda )) \zeta (\lambda )
=\sum_{s_{n_i}(\mu )=r_{m_i}(\lambda )} \overline{v_{i,q}(\mu )} \xi_i(\mu ) \zeta (\lambda )
\end{equation*}
If $\lambda \notin \mathrm{Seg}_{(n_i,n)}^n (K_i \times_c \Lambda^{m_i})$, then $(\phi_{m_i}(\langle v_{i,q},f_i \rangle )\zeta )(\lambda )=0$.
 Moreover $\{ u_{i,r} \}_r$ is a partition of unity for $\mathrm{Seg}_{(n_i,n)}^n (K_i \times_c \Lambda^{m_i})$, we obtain
\begin{eqnarray*}
\sum_{s_{n_i}(\mu )=r_{m_i}(\lambda )} \overline{v_{i,q}(\mu )} \xi_i(\mu ) \zeta (\lambda )
&=&\sum_{s_{n_i}(\mu )=r_{m_i}(\lambda )} \overline{v_{i,q}(\mu )} \xi_i(\mu ) \sum_{r} u_{i,r}^2(\lambda )\zeta (\lambda )\\
&=& \sum_{r} (\langle v_{i,q},\xi_i \rangle  \cdot u_{i,r} \cdot \langle u_{i,r},\zeta \rangle )(\lambda ).
\end{eqnarray*}
Since
$T_{n_i}(v_{i,q})^*T_{n_i}(\xi_i)=T_0(\langle v_{i,q} , \xi_i \rangle )=\psi^T_{m_i}(\phi_{m_i}
(\langle v_{i,q},\xi_i \rangle ))$, we obtain the equation
\begin{equation}
T_{n_i}(v_{i,q})^*T_{n_i}(\xi_i)=\sum_r T_{m_i}( \langle v_{i,q} , \xi_i \rangle \cdot u_{i,r} )T_{m_i}( u_{i,r} )^*. \tag{$\sharp$}
\end{equation}
For $x \in T_0(A) \cap C^*(T)$,
\begin{eqnarray*}
\lefteqn{\sigma_{\{v_{i,q} \otimes u_{i,p} \}_{p,q}}^{T,n}(x)T_{n_i}(\xi_i)}\hspace{0.5cm}\\
&=&\sum_{p,q} T_{n_i}(v_{i,q})T_{m_i}(u_{i,p})xT_{m_i}(u_{i,p})^*T_{n_i}(v_{i,q})^*T_{n_i}(\xi_i)\\
&=&\sum_{p,q,r} T_{n_i}(v_{i,q})T_{m_i}(u_{i,p})xT_{m_i}(u_{i,p})^* T_{m_i}( \langle v_{i,q} , \xi_i \rangle \cdot u_{i,r} )T_{m_i}( u_{i,r} )^* \qquad \mathrm{(by \ Eq.(\sharp ))}\\
&=&\sum_{p,q,r} T_{n_i}(v_{i,q})T_{m_i}(u_{i,p})T_{m_i}(u_{i,p})^* T_{m_i}( \langle v_{i,q} , \xi_i \rangle \cdot u_{i,r} )xT_{m_i}( u_{i,r} )^*\\
&=&\sum_{p,q,r} T_{n_i}(v_{i,q})T_{m_i}(u_{i,p} \cdot \langle u_{i,p} , \langle v_{i,q} , \xi_i \rangle \cdot u_{i,r} \rangle )xT_{m_i}( u_{i,r} )^*\\
&=&\sum_{q,r} T_{n_i}(v_{i,q})T_{m_i}( \langle v_{i,q} , \xi_i \rangle \cdot u_{i,r} )xT_{m_i}( u_{i,r} )^* \qquad \mathrm{(by \ Lemma\ \ref{lem:unity})}\\
&=&\sum_{q,r} T_{n_i}(v_{i,q} \cdot \langle v_{i,q} , \xi_i \rangle)T_{m_i}( u_{i,r} )xT_{m_i}( u_{i,r} )^*\\
&=&\sum_{r} T_{n_i}(\xi_i)T_{m_i}( u_{i,r} )xT_{m_i}( u_{i,r} )^* \qquad \mathrm{(by \ Lemma\ \ref{lem:unity})} \\
&=&T_{n_i}(\xi_i)\sigma_{\{u_{i,r}\}_r}^{T,m_i}(x).
\end{eqnarray*}
The proof is completed.
\end{proof}

For $m \in \mathbb{N}^k$, let us define an injective *-homomorphism $\pi_m :C_b(\Lambda^m ) \longrightarrow L(X_m)$ by
\begin{equation*}
(\pi_m(f)\xi )(\lambda )=f(\lambda ) \xi (\lambda ), \quad \lambda \in \Lambda^m. 
\end{equation*}
The following lemma is proved by Katsura (Lemma 1.16, 1.17 of \cite{Katsura1})

\begin{lem}
\label{lem:katsura}
\normalfont\slshape
For each $m \in \mathbb{N}$, the image $\pi_m(C_0(\Lambda^m))$ of $C_0(\Lambda^m)$ is included in the $C^*$-algebra $K(X_m)$ injectively. In particular,  
given $f \in C_c(\Lambda^m)$, there exists a orthogonal pair $\{(\xi_i, \eta_i )\}_{i=1}^L$ such that
\begin{equation*}
f=\sum_{i=1}^L \xi_i \overline{\eta_i} ,\quad \mathrm{and} \quad \pi_m(f)=\sum_{i=1}^L \theta_{\xi_i , \eta_i}
\end{equation*}
\end{lem}
 
\begin{lem}
\normalfont\slshape
 Let $\Lambda$ be a row-finite topological $k$-graph with no source.
 Set $\varphi_m^T=\psi_m^T \circ \pi_m :C_0(\Lambda^m ) \longrightarrow C^*(T)$.
 Then, for any $f \in C_0(\Lambda^m)$, the element $\varphi_m^T(f) $ belongs to the relative commutant algebra $T_0(A)' \cap C^*(T)$.
\end{lem}

\begin{proof}
For all $g \in A$, since $\phi_m(g)=\pi_m(g \circ r_m)$, 
\begin{equation*}
T_0(g)\varphi_m^T(f)
 = \psi_m^T(\phi_m(g))\psi_m^T(\pi_m(f))
 = \psi_m^T(\pi_m(f))\psi_m^T(\phi_m(g))
 = \varphi_m^T(f)T_0(g).
\end{equation*}
\end{proof}

\begin{lem}
\label{lem:turn3}
\normalfont\slshape
Let $\{ (u_i,u_i) \}_{i=1}^L$ be a orthogonal pair for degree $m$. Then for $f \in C_c(\Lambda^n)$,
\begin{equation*}
\sigma_{\{u_i\}_i}^{T,m}(\varphi_n^T(f)))=\varphi_{m+n}^T(\sum_{i=1}^L |u_i|^2 \otimes f)
\end{equation*}
\end{lem}

\begin{proof}
By Lemma \ref{lem:katsura}, there is a orthogonal pair $\{ (\xi_j, \eta_j) \}_{j=1}^M$ such that $f= \sum_{j=1}^M \xi_j \overline{\eta_j}$. Then $\{ (u_i \otimes \xi_j, u_i \otimes \eta_j) \}_{1 \le i \le L,1 \le j \le M}$ is a orthogonal pair for degree $m+n$ and
\begin{eqnarray*}
\sigma_{\{u_i\}_i}^{T,m}(\varphi_n^T(f)))
&=&\sum_{i,j}T_{m+n}(u_i \otimes \xi_j) T_{m+n}(u_i \otimes \eta_j)^* 
 = \varphi^T_{m+n}(\sum_{i,j} (u_i \otimes \xi_j)\overline{(u_i \otimes \eta_j)} ) \\
&=&\varphi^T_{m+n}(\sum_{i=1}^L |u_i|^2 \otimes f) .
\end{eqnarray*}
\end{proof}

\begin{df}
\label{df:aperiodic}
 Let $\Lambda$ be a row-finite topological $k$-graph with no source.
 We say $\Lambda$ satisfies \textit{condition (A)} if for any $v \in \Lambda^0$ and for any open neighborhood $V$ at $v$, there exist $v' \in V$ and infinite path $\alpha \in \Lambda^\infty (v')$ such that $p,q \in \mathbb{N}^k$, $p \neq q$ implies $\tau^p(\alpha) \neq \tau^q(\alpha)$. For such infinite path $\alpha$, we say that $\alpha$ is an \textit{aperiodic path}.
\end{df}

\begin{prop}
\label{prop:freeness}
\normalfont\slshape
Let $\Lambda$ be a row-finite topological $k$-graph with no source satisfying condition (A) and $T$ is injective.
For any $\epsilon >0$ and $x=\sum_{i=1}^L T_{n_{i,1}}(\xi_{i,1})T_{n_{i,2}}(\xi_{i,2})^* \in C^*(T)^{\mathrm{cpt}}$, there exist $b_1,b_2 \in C^*(T)$ and $S \in K(X_l)$ such that
\begin{equation*}
\| b_1 \| , \| b_2 \| \le 1 ,\quad b_1^*xb_2=b_1^*x_0b_2=\psi_l^T(S) ,\quad \| x_0 \| < \| b_1^* x_0 b_2 \| + \epsilon
\end{equation*}
where $x_0=\sum_{ \{1 \le i \le L | n_{i,1}=n_{i,2} \} } T_{n_{i,1}}(\xi_{i,1})T_{n_{i,2}}(\xi_{i,2})^* \in \mathcal{F}_T$.
 Moreover, if $x_0 \ge 0$ then we can take $b=b_1=b_2$, $S \ge 0$.
\end{prop}

\begin{proof}
Let us put $n=\vee_{i=1}^L (n_{i,1} \vee n_{i,2} )$. Then using row-finiteness and no source, we can suppose
\begin{equation*}
x_0=\sum_i T_n(\xi_{i,1}')T_n(\xi_{i,2}')^*.
\end{equation*}  
We remark that the sum of $x_0$ may be infinite or the components may have non-compact supports. Let us put $S_0=\sum_i \theta_{\xi_{i,1}',\xi_{i,2}'} \in K(X_n)$, then since the injectivity of $T_0$, $\| x_0 \| = \| S_0\|$ by Lemma 2.2 of \cite{KajiwaraPinzariWatatani}.
 Since $\| S_0 \|=\sup_{ \| \xi \| , \| \eta \| =1} \| \langle \xi , S_0 \eta \rangle \|_\infty$, for any $\epsilon >0,$ there exists $\xi_1 ,\xi_2 \in C_c(\Lambda^n )$ such that
\begin{equation*}
\| \xi_1 \|= \| \xi_2 \| =1 , \quad \| S_0 \| - \epsilon < \| \langle \xi_1 , S_0 \xi_2 \rangle \|_\infty
\end{equation*}
Then there is an element $v' \in \Lambda^0$ and open neighbor $V$ at $v'$ such that
\begin{equation*}
\| S_0 \| -\epsilon < | \langle \xi_1 , S_0 \xi_2 \rangle (v) | \quad (v \in V)
\end{equation*}
For $V$, using condition (A), there exists $v'$ and an aperiodic path $\alpha \in \Lambda^\infty(v')$ such that $v' = r(\alpha )$.
For each $0 \le p \le n$, $\tau^p \alpha$ is a different element of $\Lambda^\infty$, there exists $m \in \mathbb{N}^k$ such that $\alpha (p,m+p)$ is each different element of $\Lambda^m$. Let us take a disjoint open set $\{ U_p \}_{0 \le p \le n}$ such that $\alpha (p,m+p) \in U_p$.
Set $l=m+n$ and $U$ by
\begin{equation*}
U=\bigcap_{0 \le p \le n}(\mathrm{Seg}_{(p,m+p)}^l)^{-1}(U_p)=\{ \lambda \in \Lambda^l | \lambda (p,m+p) \in U_p \quad 0 \le p \le n \} .
\end{equation*}
Then $U$ is an open set and $\alpha (0,l) \in U$. Since $\Lambda^l$ is locally compact, we can suppose $U$ is relative compact. Take $Q \in C_c(U)$ such that $Q(\alpha (0,l)) =1$ and $0 \le Q \le 1$.\\
 Fix $i$ such that $n_{i,1} \neq n_{i,2}$ and set $m_{i,j}=n-n_{i,j}$.
 Let us put $K_1=\mathrm{supp}(\xi_{i,1})$, $K_{2}=\mathrm{supp}(\xi_{i,2})$, $K_3=\mathrm{supp}(\xi_1 ) \cup \mathrm{supp}(\xi_2 )$.
 By Lemma \ref{lem:construct unity}, there exist a partition of unity $\{v_{j,q} \}_q$ for $K_j$ and $\{ u_{j,p} \}$ for $\mathrm{Seg}_{(n_{i,j},n)}^n(K_j \times_c \Lambda^{m_{i,j}}) $ $(j=1,2)$ such that $\{v_{j,q} \otimes u_{j,p} \}_{p,q}$ is a partition of unity for $K_3$. By Lemma \ref{lem:turn1} and \ref{lem:turn2}, 
\begin{equation*}
\sigma_{\{v_{j,q} \otimes u_{j,p} \}_{p,q}}^{T,n}(x) T^n(\xi_j)=T^n(\xi_j) x ,
\quad
\sigma_{\{v_{j,q} \otimes u_{j,p} \}_{p,q}}^{T,n}(x) T_{n_{i,j}}(\xi_{i,j})
= T_{n_{i,j}}(\xi_{i,j}) \sigma_{\{u_{j,p} \}_{p}}^{T,m_{i,j}}(x).
\end{equation*} 
for $x \in T^0(A)' \cap C^*(T)$, $j=1,2$. On the other hand,
\begin{eqnarray*}
\lefteqn{\sigma_{\{u_{1,p} \}_{p}}^{T,m_{i,1}}(\varphi_l^T (Q)) \sigma_{\{u_{2,p} \}_{p}}^{T,m_{i,2}}(\varphi_l^T (Q))}\hspace{1cm}\\
&=&\varphi_{l+m_{i,1}}^T(\sum_{p_1} (u_{1,p_1})^2 \otimes Q)\varphi_{l+m_{i,2}}^T(\sum_{p_2} (u_{2,p_2})^2 \otimes Q) 
 \qquad \mathrm{(by \ Lemma\ \ref{lem:turn3})} \\
&=&\psi_{l+m_{i,1}}^T \Bigl( \pi_{l+m_{i,1}} (\sum_{p_1} (u_{1,p_1})^2 \otimes Q ) \Bigr) \psi_{l+m_{i,2}}^T \Bigl( \pi_{l+m_{i,2}} (\sum_{p_2} (u_{2,p_2})^2 \otimes Q ) \Bigr)\\
&=&\psi_{(l+m_{i,1}) \vee (l+m_{i,2})}^T\Big( \pi_{l+m_{i,1}} (\sum_{p_1} (u_{1,p_1})^2 \otimes Q ) \otimes 1)
\pi_{l+m_{i,2}} (\sum_{p_2} (u_{2,p_2})^2 \otimes Q ) \otimes 1\Big).
\end{eqnarray*}
In the last equation, we used Nica covariance. If we suppose $\mu (m_{i,1}, l+m_{i,1}) \in U$ and $\mu (m_{i,2} , l+m_{i,2}) \in U$ for $\mu \in \Lambda^{(l+m_{i,1}) \vee (l+m_{i,2})}$, then this contradicts  $U_{(m_{i,1} \vee m_{i,2}) - m_{i,1}} \cap U_{(m_{i,1} \vee m_{i,2}) - m_{i,2}} = \emptyset$. Hence we obtain
\begin{eqnarray*}
\lefteqn{\varphi_l^T(Q)T_n(\xi_1 )^* T_{n_{i,1}}(\xi_{i,1})T_{n_{i,2}}(\xi_{i,2})^* T_n(\xi_2 ) \varphi_l^T(Q)}\hspace{1cm}\\
&=&T_n(\xi_1 )^* T_{n_{i,1}}(\xi_{i,1}) \Bigl( \sigma_{\{u_{1,p} \}_{p}}^{T,m_{i,1}}(\varphi_l^T (Q)) \sigma_{\{u_{2,p} \}_{p}}^{T,m_{i,2}}(\varphi_l^T (Q)) \Bigr) T_{n_{i,2}}(\xi_{i,2})^* T_n(\xi_2 )\\
&=&0
\end{eqnarray*}
for $n_{i,1} \neq n_{i,2}$. Set $b_j=T_n(\xi_j)\varphi_l^T(Q)$ $(j=1,2)$. Then we get $b_1^*x_0 b_2=b_1^* x b_2$ and $\| b_1 \| , \| b_2 \| \le 1$.
 Finally we have
\begin{eqnarray*}
\| x_0 \| - \epsilon 
& < & | \langle \xi_1 , S_0  \xi_2 \rangle (v' ) | 
=|Q(\alpha (0,l) ) \langle \xi_1 , S_0 \xi_2 \rangle (v' ) Q(\alpha (0,l) ) | \\
 &\le&  \| Q (\langle \xi_1 , S_0 \xi_2 \rangle \circ r_l) Q \| 
= \| b_1^* x_0 b_2 \| .
\end{eqnarray*}
The proof is completed.
\end{proof}

\begin{thm}
\label{thm:CK-unique}
\normalfont\slshape
Let $\Lambda$ be a row-finite topological $k$-graph with no source. If $\Lambda$ satisfies Condition (A) and $T$ is an injective Cuntz-Pimsner covariance representation, then the natural map $\mathcal{O}(\Lambda ) \longrightarrow C^*(T)$ is an isomorphism.
\end{thm}

\begin{proof}
Take $x \in C^*(T)^{\mathrm{cpt}}$. Let $x,x_0$ be represent as in Proposition \ref{prop:freeness}. For the proof of theorem, we enough to show $\| x_0 \| \le \| x \|$ since $C^*(T)^{\mathrm{cpt}}$ is dense in $C^*(T)$. But this inequation follows from Proposition \ref{prop:freeness}.
\end{proof}

The following proposition will use in Lemma \ref{lem:connect} (see also Proposition 5.10 of \cite{Katsura1}).
\begin{prop}
\label{prop:katsura5.10}
\normalfont\slshape
Let $\Lambda$ be a row-finite topological $k$-graph and satisfies condition (A) and $T$ is injective.
For $\epsilon>0,\ x \in C^*(T)^{\mathrm{cpt}}$, we define $x_0 \in \mathcal{F}_T$ as in Proposition \ref{prop:freeness}. If $x_0\ge0$, then there exist $a \in C^*(T)$ and a positive function $f \in C_0(\Lambda^0)$ such that
\begin{equation*}
\| f \| = \| x_0 \| ,\quad \| a^*xa-T_0(f) \| < \epsilon .
\end{equation*}
\end{prop}

\begin{proof}
By Proposition \ref{prop:freeness}, there exist $b \in C^*(T)$ and $0 \le S \in K(X_l)$ such that 
\begin{equation*}
\| b \| \le 1 ,\quad b^* x b=b^* x_0 b=\psi_l^T(S) ,\quad \| x_0 \| < \| b^* x_0 b \| + \frac{\epsilon}{2}
\end{equation*}
 Then there is $\zeta \in X_l$ such that $\| \zeta \| =1$ and $\| \langle \zeta , S \zeta \rangle \| > \| S \| -\epsilon /2$. 
 Set $g=\langle \zeta , S \zeta \rangle \in C_0(\Lambda^0)$ and $a=bT_l(\zeta )$.
 Then $T_0(g)=a^*x_0a$. Since $\| g \| \le \| S\| \le \| x_0 \|$ and 
\begin{equation*}
\| x_0 \| < \| b^*x_0b \| + \frac{\epsilon}{2} = \| S\| + \frac{\epsilon}{2} < \| \langle \zeta , S \zeta \rangle \| + \epsilon = \| g \| + \epsilon.
\end{equation*}
Therefore we get $| \| g \| - \| x_0 \| | <\epsilon /2$. Hence $f=\|x_0\| g/\| g\|$ satisfies that $\| f \| = \| x_0 \|$ and 
\begin{eqnarray*}
\| a^*xa-T_0(f) \| &\le& \| a^*x_0a - T_0(g) \| + \| T_0(g)-T_0(f) \| 
< | \| g \| - \| x_0 \| | < \epsilon 
\end{eqnarray*}
\end{proof}


\section{Simplicity and Purely infiniteness}
\label{sec:simple}
In this section, we give criteria when $\mathcal{O}(\Lambda )$ is simple and purely infinite. Our terms and analysis are based on the works of Katsura (\cite{Katsura3}, \cite{Katsura4}), Kajiwara-Pinzari-Watatani(\cite{KajiwaraPinzariWatatani}), and Kajiwara-Watatani (\cite{KajiwaraWatatani}).

\begin{df}
Let $\Lambda$ be a row-finite topological $k$-graph with no source and $X$ be a product system associated with $\Lambda$.
\begin{enumerate}
 \item{For an ideal $I$ of $C_0(\Lambda^0)$, $I$ is \textit{$X$-invariant} if for any $f \in C_0(\Lambda^0)$ and $\xi ,\eta \in X_m$, $\langle \xi , f \cdot \eta \rangle$ belongs to $I$. $I$ is \textit{$X$-saturated} if there is $1 \le i \le k$ such that $\langle \xi , f \cdot \eta \rangle \in I$ for any $\xi ,\eta \in X_{e_i}$, then $f \in I$.}
 \item{Let $\Omega$ be a subset of $\Lambda^0$. We say that $\Omega$ is \textit{positively invariant} if for $\lambda \in \Lambda$, $s(\lambda )\in \Omega$ implies $r(\lambda ) \in \Omega$, and $\Omega$ is \textit{negatively invariant} if for any $v \in \Omega$ and $1 \le i \le k$, there exists $\lambda_i \in \Lambda^{e_i}$ such that $v=r(\lambda_i)$ and $s(\lambda_i) \in \Omega$. We say that $\Omega$ is \textit{invariant} if $\Omega$ is positively and negatively invariant.}
 \item{$\Lambda$ is said to be \textit{minimal} if there exist no closed invariant set other that $\emptyset$ or $\Lambda^0$.}
\end{enumerate}
\end{df}

Let $I$ be an ideal of $C_0(\Lambda^0)$. Then there exists a closed subset $\Omega_I$ of $\Lambda^0$ such that $I=C_0(\Lambda^0 \setminus \Omega_I )$.

\begin{prop}
\label{prop:inv-sat}
\normalfont\slshape
 Under the above assumption,
\begin{enumerate}
 \item{$I$ is $X$-invariant if and only if $\Omega_I$ is positively invariant.}
 \item{$I$ is $X$-saturated if and only if $\Omega_I$ is negatively invariant.}
\end{enumerate}
\end{prop}

\begin{proof}
Suppose $I$ is $X$-invariant and $s_m(\lambda ) \in \Omega_I$. Let $0 \le a \in I$, then $\langle \xi , f \cdot \xi \rangle \in I$ for any $\xi \in X_m$. Take $\xi \in X_m$ such that $\xi(\lambda )=1$. Then $r_m(\lambda ) \in \Omega_I$ since
\begin{equation*}
0= \langle \xi , f \cdot \xi \rangle (s_m(\lambda ))=\sum_{s_m(\mu )=s_m(\lambda )} f(r_m(\mu ))|\xi (\mu )|^2 \ge f(r_m(\lambda ))|\xi (\lambda )|^2=f(r_m(\lambda )) \ge 0.
\end{equation*}
Conversely, we suppose $\Omega_I$ is positively invariant. For each $f \in I$ and $\lambda \in \Lambda^m$ such that $s_m(\lambda ) \in \Omega_I$,
\begin{equation*}
\langle \xi, f \cdot \eta \rangle (s_m(\lambda ))=\sum_{s_m(\mu )=s_m(\lambda )}\overline{\xi (\mu )} f(r_m(\mu )) \eta (\mu )=0
\end{equation*}
since $r_m(\mu ) \in \Omega_I$ for all $\mu \in \Lambda^m$ such that $s_m(\mu )=s_m(\lambda )$. Hence we have done the proof of (1).

Next, we shall prove (2).
 First we suppose $\Omega_I$ is negative invariant.
 Let us suppose that there is $1 \le i \le k$ such that $\langle \xi, f \cdot \eta \rangle \in I$ for all $\xi ,\eta \in X_{e_i}$.
 For any $v \in \Omega_I$, there is $\lambda \in \Lambda^{e_i}$ such that $v=r_{e_i}(\lambda )$ and $s_{e_i}(\lambda ) \in \Omega_I$.
 Let us take $\xi ,\eta \in X_{e_i}$ such that $f(\lambda )=1 = g(\lambda )$ and $f(\mu )=0 =g(\mu )$ for $s_{e_i}(\lambda )=s_{e_i}(\mu )$ and $\lambda \neq \mu$. Then
\begin{equation*}
0= \langle \xi ,f \cdot \eta \rangle (v)= \sum_{s_{e_i}(\mu)=v} \overline{\xi (\mu )}f(r_{e_i}(\mu )) \eta(\mu )=f(r_{e_i}(\lambda ))=f(v).
\end{equation*}
This implies $f \in I$.
 Next we shall prove the "only if" part.
 We suppose that there exist $v_0 \in \Omega_I$ and $1 \le i \le k$ such that $s_{e_i}(\lambda ) \notin \Omega_I$ for any $\lambda \in \Lambda^{e_i}(v_0)$ and we shall induce a contradiction.
 Then for any $\lambda \in \Lambda^{e_i}(v_0)$, there is a open neighborhood $V_\lambda$ of $s_{e_i}(\lambda )$ in $\Lambda^0$ such that $V_\lambda \subset \Lambda^0 \setminus \Omega_I$.
 Set $V=\cup_{\lambda \in \Lambda^{e_i}(v_0)} V_\lambda$.
 Then there exists an open neighborhood $W$ of $v_0$ such that $s_{e_i}(r_{e_i}^{-1}(W)) \subset V$ by Lemma 1.21 of \cite{Katsura1}.
 Define a function $f \in A$ such that $f \ge 0,\ \mathrm{supp}(f) \subset W$ and $f(v_0)=1$.  
 We suppose that there exist $v \in \Lambda^0$ and $\xi ,\eta \in X_{e_i}$ such that $\langle \xi ,f \cdot \eta \rangle (v) \neq 0$.
 Then there is $\mu \in \Lambda^{e_i}$ such that $v=s_{e_i}(\mu )$ and $f(r_{e_i}(\mu )) \neq 0$.
 Thus $r_{e_i}(\mu ) \in W$. From a choice of $W$, $v \notin \Omega_I$ holds.
 Hence we conclude $\langle \xi , f \cdot \eta \rangle \in I$.
 The assumption of $X$-saturated implies $f \in I$, however $f(v_0)=1$.
 This is a contradiction.
 
 Hence we completed the proof.
\end{proof}

Let $J$ be an ideal of $\mathcal{O}(\Lambda )$. Since $t_0$ is injective, for $J_0:=J \cap t_0(C_0(\Lambda^0))$ there exists a closed subset $\Omega_{I_0}$ of $\Lambda^0$ such that $J_0=t_0(C_0(\Lambda^0 \setminus \Omega_{I_0}))$. Let us put $I_0=C_0(\Lambda^0 \setminus \Omega_{I_0} )$ of an ideal of $C_0(\Lambda^0)$.

\begin{lem}
\normalfont\slshape
$I_0$ is $X$-invariant and $X$-saturated. Hence $\Omega_{I_0}$ is a closed invariant subset of $\Lambda^0$.
\end{lem} 

\begin{proof}
Since $t_0(f) \in J_0$ for $f \in I_0$, $t_0(\langle \xi ,f \cdot \eta \rangle ) = t_m(\xi )^*t_0(f)t_m(\eta ) \in J_0$ for $\xi ,\eta \in X_m$.
 Hence $\langle \xi ,f \cdot \eta \rangle \in I_0$ and this implies that $I_0$ is $X$-invariant.\\
Let us suppose that there exists $1 \le i \le k$ such that $\langle \xi , f \cdot \eta \rangle \in I_0$ for $\xi ,\eta \in X_{e_i}$.
 Then for any $S_1,S_2 \in K(X_{e_i})$, $\psi_{e_i}(S_1)t_0(f)\psi_{e_i}(S_2) \in J_0$.
 Since $\Lambda$ is row-finite without source, $t_0(f) \in \psi_{e_i}(K(X_{e_i}))$ and using an approximate unit of $K(X_{e_i})$, we can conclude $t_0(f) \in J_0$, hence $f \in I_0$.
\end{proof}

\begin{lem}
\label{lem:invideal}
Let $\Omega$ be a non-trivial closed invariant subset of $\Lambda^0$. 
 Set $I_0=C_0(\Lambda^0 \setminus \Omega )$ and $J$ be an ideal generated by $t_0(I_0)$ in $\mathcal{O}(\Lambda )$.
 Then $J \cap t_0(A)=t_0(I_0)$ holds.
\end{lem}

\begin{proof}
Let us put
\begin{equation*}
J_{\mathrm{alg}}=\mathrm{span} \{ t_{m_1}(\xi_1)t_{m_2}(\xi_2)^* t_0(f) t_{n_1}(\eta_1)t_{n_2}(\eta_2)^* | f \in I_0,\ 
\xi_i \in X_{m_i} ,\ \eta_i \in X_{n_i} \}
\end{equation*}
which is dense in $J$. Then we shall show $J_{\mathrm{alg}} \cap t_0(A) \subset t_0(I_0)$.
\begin{equation*}
x_0 = \sum_{j=1}^L t_{m_{j,1}}(\xi_{j,1})t_{m_{j,2}}(\xi_{j,2})^* t_0(f_j) t_{n_{j,1}}(\eta_{j,1})t_{n_{j,2}}(\eta_{j,2})^* \in J_{\mathrm{alg}} \cap t_0(A).
\end{equation*}
By using the faithful conditional expectation $\Psi$ onto $\mathcal{O}(\Lambda )^\gamma$, we can assume $m_{j,1}+ n_{j,1}=m_{j,2}+n_{j,2}$ for every $1 \le j \le L$.
 Then there is a large $l \in \mathbb{N}^k$ such that for any $\xi, \eta \in X_l$,
\begin{equation*}
t_l(\xi )^* x_0 t_l(\eta )=\sum_{j=1}^L t_{l_j}(\xi_j')^* t_0(f_j) t_{l_j}(\eta_j')
=t_0(\sum_{j=1}^L \langle \xi_j' , f_j \cdot \eta_j' \rangle ).
\end{equation*}
where $l_j \in \mathbb{N}^k$ are suitable elements and $\xi_j', \eta_j' \in X_{l_j}$.
 We remark $\sum_{j=1}^L \langle \xi_j' , f_j \cdot \eta_j' \rangle \in I_0$ since $I_0$ is $X$-invariant by Proposition \ref{prop:inv-sat}.
 On the other hand, there exists $f_0 \in A$ such that $x_0=t_0(f_0)$.
 Since $I_0$ is $X$-saturated, we obtain $x_0=t_0(f_0) \in t_0(I_0)$.

Next, we prove $J \cap t_0(A) = t_0(I_0)$. The inclusion $J \cap t_0(A) \supset t_0(I_0)$ is obvious.
Take $x \in J \cap t_0(A)$. For any $\epsilon >0$, there is 
\begin{equation*}
x_0=\sum_{j=1}^L t_{m_{j,1}}(\xi_{j,1})t_{m_{j,2}}(\xi_{j,2})^* t_0(f_j) t_{n_{j,1}}(\eta_{j,1})t_{n_{j,2}}(\eta_{j,2})^* \in J_{\mathrm{alg}}
\end{equation*}
such that $\| x-x_0 \| < \epsilon$.
 By using the faithful conditional expectation $\Psi$ again, we can assume $m_{j,1}+ n_{j,1}=m_{j,2}+n_{j,2}$ for every $1 \le j \le L$.
 For a large $l \in \mathbb{N}^k$ and any $\xi , \eta \in X_l$ such that $\| \xi \|$, $\| \eta \| \le 1$, we obtain $\| t_l(\xi )^*x t_l(\eta ) - t_l(\xi )^*x_0 t_l(\eta ) \| < \epsilon$ and $t_l(\xi )^* x_0 t_l(\eta ) \in J_{\mathrm{alg}} \cap t_0(A) \subset t_0(I_0)$.
 Hence we obtain $t_l(\xi )^*x t_l(\eta ) \in \overline{t_0(I_0)}=t_0(I_0)$. Since $x \in t_0(A)$ and $I_0$ is $X$-saturated, we obtain $x \in t_0(I_0)$. Therefore we finished.
\end{proof}

Next, we introduce an orbit space followed by Katsura (\cite{Katsura3}). 
\begin{df}
 Let $\Lambda$ be a row-finite topological $k$-graph with no source.
 For $v \in \Lambda^0$, define a subset $\mathrm{Orb}^+(v)$ of $\Lambda^0$ by
\begin{equation*}
\mathrm{Orb}^+(v)=\bigcup_{n \in \mathbb{N}^k} r_n( s_n^{-1}(v)). 
\end{equation*}
 For $v \in \Lambda^0$ and $\alpha \in \Lambda^\infty (v)$, define subsets $\mathrm{Orb}^{-}(v,\alpha )$ and $\mathrm{Orb}(v,\alpha )$ of $\Lambda^0$ by
\begin{equation*}
\mathrm{Orb}^{-}(v,\alpha )=\{ \alpha(m) | m \in \mathbb{N}^k \} ,\quad 
\mathrm{Orb}(v,\alpha )=\bigcup_{v' \in \mathrm{Orb}^{-}(v,\alpha )} \mathrm{Orb}^+(v'). 
\end{equation*}
The subsets $\mathrm{Orb}^+(v),\ \mathrm{Orb}^{-}(v,\alpha ),\ \mathrm{Orb}(v,\alpha )$ are said to be the \textit{positive orbit space} of $v \in \Lambda^0$, \textit{negative orbit space} of $(v,\alpha )$, \textit{orbit space} of $(v,\alpha )$, respectively.
\end{df}

\begin{lem}
\label{lem:invariant}
\normalfont\slshape
Let $\Lambda$ be a row-finite topological $k$-graph with no source. 
\begin{enumerate}
\item{Let $\Omega$ be a closed invariant subset of $\Lambda^0$ and take $v \in \Omega$. Then there is $\alpha \in \Lambda^\infty(v)$ such that $\mathrm{Orb}(v,\alpha ) \subset \Omega$.}
\item{For any $v \in \Lambda^0$ and $\alpha \in \Lambda^\infty(v)$, $\mathrm{Orb}(v,\alpha ) \subset \Lambda^0$ is an invariant subset.}
\end{enumerate}
\end{lem} 

\begin{proof}
(1) First, we shall show the case $k<\infty $. We remark that we need not to assume $\Omega$ is closed in this case. Let $\Omega$ be an invariant subset of $\Lambda^0$. Fix $v \in \Omega$. By negativity and positivity of $\Omega$ and factorization property of $\Lambda$, we can easily construct $\alpha \in \Lambda^\infty$ such that $\alpha (m) \in \Omega$ for any $m \in \mathbb{N}^k$. Using the positivity of $\Omega$, $\mathrm{Orb}(v,\alpha ) \subset \Omega$.

Next we shall show the case $k=\infty$. For $m \in \mathbb{N}^k$ and an open set $U \subset \Lambda^m$, let us define $Z(U)=\{ \alpha \in \Lambda^\infty | \alpha (0,m) \in U \}$.
 We shall define a topology of $\Lambda^\infty$ such that $\{ Z(U)\ |\ m \in \mathbb{N}^k,\ U $ is open subset of $\Lambda^m \}$ forms the open basis.
 For $v \in \Omega$, $\Lambda^\infty(v)$ is a topological space by the relative topology of $\Lambda^\infty$.
 On the other hand, for $n \le m$, we shall define $\pi_{n,m} : \Lambda^m(v) \longrightarrow \Lambda^n(v)$ by $\pi_{n,m}=\mathrm{Seg}_{(0,n)}^m$.
 Then $\{ \Lambda^m(v), \pi_{n,m} \}_{m \in \mathbb{N}^\infty }$ is an inverse system.
 Define a map $F: \Lambda^\infty (v) \longrightarrow \varprojlim_{m \in \mathbb{N}^\infty} \Lambda^m(v)$ by $F(\alpha )=\{ \alpha (0,m) \}_{m \in \mathbb{N}^\infty}$, then this map induces a homeomorphism.
 Hence $\Lambda^\infty (v)$ is compact (Hausdorff) space since $\Lambda$ is row-finite.
 For each $1 \le p < \infty$, set $\mathbb{N}^\infty_p=\{ m \in \mathbb{N}^\infty \ |\ m_{(i)}=0$ for $i > p \}$ and $E_p=\{ \alpha \in \Lambda^\infty (v)\ |\ \alpha (m) \in \Omega$ for $ m \in \mathbb{N}_p^\infty \}$.
 Then $E_p$ is non-empty closed subset of $\Lambda^\infty(v)$ since the first part of this proof and $\Omega$ is closed. Since $\Lambda^\infty(v)$ is compact and $E_1 \supset E_2 \supset \cdots $, the union $\cap_{p=1}^\infty E_p$ is not empty and $\alpha \in \cap_{p=1}^\infty E_p$ is just what we want.

(2) Take $\lambda \in \Lambda$ such that $s(\lambda ) \in \mathrm{Orb}(v,\alpha )$. Then there is $v' \in \mathrm{Orb}^-(v,\alpha )$ such that $s(\lambda ) \in \mathrm{Orb}^+(v')$. Thus there is $\mu \in \Lambda$ such that $s(\lambda )=r( \mu )$ and $s(\mu )=v'$. Then $\lambda \mu \in \mathrm{Orb}^+(v')$ and this implies $r(\lambda )=r(\lambda \mu ) \in \mathrm{Orb}^+(v') \subset \mathrm{Orb}(v,\alpha )$. Hence $\mathrm{Orb}(v,\alpha )$ is positively invariant. Take $w \in \mathrm{Orb}(v,\alpha )$. Then there exist $v' \in \mathrm{Orb}^-(v, \alpha )$ and $\lambda \in \mathrm{Orb}^+(v')$ such that $w=r(\lambda),\ v'=s(\lambda )$. Fix $1 \le i \le k$. First, we suppose $d(\lambda )_{(i)} >0$. Then $\lambda(0,e_i)$ satisfies $s(\lambda (0,e_i)) \in \mathrm{Orb}(v, \alpha )$. If $d(\lambda )_{(i)}=0$, then $\mu= \lambda \alpha (m,m+e_i) \in \mathrm{Orb}(v,\alpha)$ where $v'=\alpha (m)$. Then $\mu(0,e_i) \in \Lambda^{e_i}$ satisfies $w=r(\mu (0,e_i))$ and $s(\mu (0,e_i)) \in \mathrm{Orb}(v,\alpha )$. Hence we conclude that $\mathrm{Orb}(v,\alpha )$ is negatively invariant.
\end{proof}

\begin{thm}
\normalfont\slshape
Let $\Lambda$ be a row-finite topological $k$-graph with no source satisfying condition (A). Then
the following conditions are equivalent:
\begin{enumerate}
 \item[(i)]{The $C^*$-algebra $\mathcal{O}(\Lambda )$ is simple}
 \item[(ii)]{$\Lambda$ is minimal.}
 \item[(iii)]{The orbit space $\mathrm{Orb}(v,\alpha )$ is dense in $\Lambda^0$ for every $v \in \Lambda^0$ and $\alpha \in \Lambda^\infty(v)$.} 
\end{enumerate}
\end{thm}

\begin{proof}
(i)$\iff$(ii): Let us put $\pi : \mathcal{O}(\Lambda ) \longrightarrow \mathcal{O}(\Lambda )/J$ be a quotient map and $T_0=\pi \circ t_0$. First, we discuss the necessary condition. Then $J_0=C_0(\Lambda^0)$ or $J_0=\{ 0 \}$ by an assumption. If $J_0=C_0(\Lambda^0)$, then $T_0=0$ and this implies $J=\mathcal{O}(\Lambda )$.
If $J_0=\{ 0 \}$, then $\pi$ is identity map on $T_0(C_0(\Lambda^0))$ and $T_0$ is faithful. By condition (A) and a standard argument using Theorem \ref{thm:CK-unique}, $\pi$ is injective and this implies $J=\{ 0\}$. 
 Conversely if $\Lambda$ is not minimal, then there exists a non-trivial closed invariant subset $\Omega$ of $\Lambda^0$.
 Let us put $J \neq \{ 0 \}$ as an ideal of $\mathcal{O}(\Lambda )$ generated by $t_0(C_0(\Lambda^0 \setminus \Omega ))$.
 From Lemma \ref{lem:invideal},
\begin{equation*}
\mathcal{O}(\Lambda )/J \supset t_0(A)+ J /J \cong t_0(A)/t_0(A) \cap J \cong t_0(C_0(\Omega )) \neq \{ 0 \}.
\end{equation*}
 This says $J \neq \mathcal{O}(\Lambda )$, hence $\mathcal{O}(\Lambda )$ is not simple.\\
(ii)$\iff$(iii): This follows from Lemma \ref{lem:invariant} and if $\Omega$ is invariant then the closure $\overline{\Omega}$ is an invariant closed subset.
\end{proof}

Next, we discuss the purely infiniteness of $\mathcal{O}(\Lambda )$. We say that a simple $C^*$-algebra is \textit{purely infinite} if every non-zero hereditary $C^*$-subalgebra has an infinite projection.

\begin{lem}
\label{lem:connect}
\normalfont\slshape
Let $\Lambda$ be a row-finite topological $k$-graph with no source satisfying condition (A) and $v_0 \in \Lambda^0$ be an element of $\Lambda^0$ with $\overline{\mathrm{Orb}^+(v_0)}=\Lambda^0$. For a non-zero positive element $x \in \mathcal{O}(\Lambda )$, there exist $a \in \mathcal{O}(\Lambda )$ and $f \in C_0(\Lambda^0)$ which is 1 on some neighborhood $V_0$ of $v_0$ such that $\| a^*xa - t_0(f) \| < 1/2$.
\end{lem}

\begin{proof}
This statement is similar to Lemma 1.12 of \cite{Katsura4} and Katsura's proof works by some modifications (use Proposition \ref{prop:katsura5.10}).
\end{proof}

\begin{df}
 Let $\Lambda$ be a row-finite topological $k$-graph with no source.
 For $n,m \in \mathbb{N}^k$ and two subsets $U \subset \Lambda^n$ and $U' \subset \Lambda^m$, we define $U \pitchfork U' \subset \Lambda^{n \wedge m}$ by $U \pitchfork U' = \mathrm{Seg}^n_{(0,n \wedge m)}(U) \cap \mathrm{Seg}^m_{(0,n \wedge m)}(U')$ where $n \wedge m \in \mathbb{N}^k$ is defined by $(n \wedge m)_{(i)}=\min \{ n_{(i)}, m_{(i)} \}$.
\end{df}

\begin{lem}
\label{lem:perp}
\normalfont\slshape
 Let $\Lambda$ be a row-finite topological $k$-graph with no source.
Let $U \subset \Lambda^n$ and $U' \subset \Lambda^m$ be open sets satisfying $U \pitchfork U' = \emptyset$. Then for any $\xi \in C_c(U) \subset X_n$ and $\eta \in C_c(U') \subset X_m$, we have $t_n(\xi )^* t_m(\eta )=0$
\end{lem}

\begin{proof}
Take an approximate unit $\{ h_i \}_{i \in \mathcal{A}}$ in $C_0(\Lambda^0)$. Then, for $i,j \in \mathcal{A}$, 
\begin{eqnarray*}
t_n(\xi \cdot h_i )^* t_m(\eta \cdot h_j) 
&=& (t_n(\xi ) t_0(h_i ))^* (t_m(\eta ) t_0 (h_j ))\\
&=& (t_n(\xi ) \psi_{(n \vee m) - n}(\phi_{(n \vee m) - n} (h_i )))^* (t_m(\eta ) \psi_{(n \vee m) - m }(\phi_{(n \vee m) - m } (h_j ))).
\end{eqnarray*}
For $\xi ' \in X_{(n \vee m) - n} , \ \eta ' \in X_{(n \vee m) - m}$, $t_{n \vee m }(\xi \otimes \xi ')^* t_{n \vee m }(\eta \otimes \eta ' )=0$ implies $t_n(\xi \cdot h_i )^* t_m(\eta \cdot h_j )=0$ for any index $i,j \in \mathcal{A}$. Hence $t_n(\xi)^* t_m(\eta )=0$
\end{proof}

\begin{df}
 Let $\Lambda$ be a row-finite topological $k$-graph with no source.
 We say that a non-empty open subset $V$ of $\Lambda^0$ is a \textit{contracting open set} if its closure $\overline{V}$ is compact and there exist finitely many non-empty open subsets $U_i \subset \Lambda^{n_i}$ for $i=0,1,\cdots ,m$ with $n_i \in \mathbb{N}^k \setminus \{ 0 \}$ satisfying
\begin{enumerate}
 \item{$r_{n_i}(U_i)\subset V$ for $i=0,1,\cdots ,m$}
 \item{$U_i \pitchfork U_j= \emptyset$ for $i,j \in \{0,1,\cdots ,m\}$ with $i \neq j$}
 \item{$\overline{V} \subset \bigcup_{i=1}^m s_{n_i}(U_i)$} 
\end{enumerate}
\end{df}

\begin{df}
We say that a row-finite topological $k$-graph $\Lambda$ with no source is \textit{contracting at} $v_0 \in \Lambda^0$ if $\overline{\mathrm{Orb}^+(v_0)}=\Lambda^0$ and any neighborhood $V_0$ of $v_0$ contains a contracting open set $V \subset V_0$. We simply say that $\Lambda$ is \textit{contracting} if $\Lambda$ is contracting at some $v_0 \in \Lambda^0$.
\end{df}

\begin{thm}
\label{thm:purelyinf}
\normalfont\slshape
 Let $\Lambda$ be a row-finite topological $k$-graph with no source satisfying condition (A).
 If $\Lambda$ is minimal and contracting, then the $C^*$-algebra $\mathcal{O}(\Lambda )$ is simple and purely infinite.
\end{thm}

\begin{proof}
This proof is same as Theorem A of \cite{Katsura4}.
\end{proof}

\section{Examples}
\label{sec:Examples}

 In this section, we construct a new topological $k$-graph from a higher-rank graph and covering maps on $\mathbb{T}$ whose $C^*$-algebras include Cuntz's $ax+b$-semigroup $C^*$-algebra over $\mathbb{N}$.
 
 First, we recall Katsura's construction from the directed graph $E=(E^0,E^1,s,r)$ and two maps $m:E^1 \longrightarrow \mathbb{Z}$ and $n : E^1 \longrightarrow \mathbb{Z}_+=\mathbb{N} \setminus \{ 0 \}$.
 We define two continuous maps $\widetilde{s},\widetilde{r} : E^1 \times \mathbb{T} \longrightarrow E^0 \times \mathbb{T}$ by $\widetilde{s}(e,z)=(s(e), z^{n(e)})$ and $\widetilde{r}(e,z)=(r(e), z^{m(e)})$. Then $E \times_{n,m} \mathbb{T} = (E^0 \times \mathbb{T} , E^1 \times \mathbb{T} , \widetilde{s} , \widetilde{r} )$ is a topological graph in the sense of Katsura. We shall consider a higher-rank version of $E \times_{n,m} \mathbb{T}$.

By Theorem 2.1 and Theorem 2.2 of \cite{FowlerSims}, a topological $k$-graph $\Lambda$ is charactrized by $\Lambda^{e_i}$ $(1 \le i \le k)$ and homeomorphisms $T_{i,j} : \Lambda^{e_i} \times_c \Lambda^{e_j} \longrightarrow \Lambda^{e_j} \times_c \Lambda^{e_i}$ $(1 \le i < j \le k)$ preserving range and source maps respectively and have the \textit{hexagonal condition}
\begin{equation*}
(T_{j,l} \otimes 1_i)(1_j \otimes T_{i,l})(T_{i,j} \otimes 1_l)
=(1_l \otimes T_{i,j})(T_{i,l} \otimes 1_j)(1_i\otimes T_{j,l})
\end{equation*}
on $\Lambda^{e_i} \times_c \Lambda^{e_j} \times_c \Lambda^{e_l}$ for $k \ge 3$ and each $1 \le i < j < l \le k$.
   
\begin{prop}
\label{prop:newgraph}
\normalfont\slshape
Let $\Gamma$ be a (discrete) row-finite higher-rank graph with no source defined by the maps $T_{i,j}: \Gamma^{e_i} \times_c \Gamma^{e_j} \longrightarrow \Gamma^{e_j} \times_c \Gamma^{e_i}$, $n_i : \Gamma^{e_i} \longrightarrow \mathbb{Z}_+$ and $m_i : \Gamma^{e_i} \longrightarrow \mathbb{Z} \setminus \{ 0 \}$.
 For each $(\lambda_1, \lambda_2) \in \Gamma^{e_i} \times_c \Gamma^{e_j},\ (\lambda_1^{ij} , \lambda_2^{ij}) \in \Gamma^{e_j} \times_c \Gamma^{e_i}$ such that $T_{i,j}(\lambda_1 , \lambda_2 )=( \lambda_1^{ij} , \lambda_2^{ij} )$, if $n_i, m_i$ satisfy the following relations (i),(ii), then there exists a homomorphism
\begin{equation*}
\widetilde{T}_{i,j} : \Gamma^{e_i} \times_{n_i, m_i} \mathbb{T} \longrightarrow \Gamma^{e_j} \times_{n_j, m_j} \mathbb{T}
\end{equation*}
such that this map induces a row-fintie topological $k$-graph $\Lambda_{(\Gamma ,m ,n)}$ with no source;
\begin{enumerate}
 \item[(i)]{ $\mathrm{gcd}(|m_i(\lambda_1)m_j(\lambda_2)|, n_i(\lambda_1)n_j(\lambda_2))=1$}
 \item[(ii)]{ $m_i(\lambda_1)m_j(\lambda_2)=m_j(\lambda_1^{ij}) m_i(\lambda_2^{ij})$,  $n_i(\lambda_1)n_j(\lambda_2)=n_j(\lambda_1^{ij}) n_i(\lambda_2^{ij})$}
\end{enumerate}
where $\mathrm{gcd}(m,n)$ is the greatest common divisor for positive integers $m,n$.
\end{prop}

\begin{proof}
Given $p \ge 1$ and $\lambda_{j} \in \Gamma^{e_{i_j}}$ $(1 \le j \le p)$, let us define 
\begin{equation*}
\mathbb{T}_{(\lambda_{1},\cdots , \lambda_{p}) }
=\{ (z_1 ,\cdots ,z_p ) \in \mathbb{T}^p | z_j^{n_{i_j}(\lambda_j)}=z_{j+1}^{m_{i_{j+1}}(\lambda_{j+1}) } \ (1 \le j \le p-1) \}
\end{equation*}
and define a range map $r_{(\lambda_{1}, \cdots \lambda_{p})}$ and a source map $s_{(\lambda_{1}, \cdots \lambda_{p})}$ by
\begin{equation*}
r_{(\lambda_{1}, \cdots \lambda_{p})}(z_1,\cdots ,z_p)=z_1^{m_{i_1}(\lambda_{1})}, \quad s_{(\lambda_{1}, \cdots \lambda_{p})}(z_1,\cdots ,z_p)=z_p^{n_{i_p}(\lambda_{p})}
\end{equation*}

For the fixed $(\lambda_i, \lambda_j) \in \Gamma^{e_i} \times_c \Gamma^{e_j}$ with $T^{e_i,e_j}(\lambda_1, \lambda_2)=(\lambda_1^{ij}, \lambda_2^{ij} )$, we enough to construct a homeomorphism $S_{i,j} : \mathbb{T}_{(\lambda_1, \lambda_2)} \longrightarrow \mathbb{T}_{(\lambda_1^{ij}, \lambda_2^{ij})}$ preserving range and source maps respectively.
 Fix $z \in \mathbb{T}$ and consider the elements $(z_1,z_2) \in \mathbb{T}_{(\lambda_1 , \lambda_2 )}$ and $(w_1,w_2) \in \mathbb{T}_{(\lambda_1^{ij} , \lambda_2^{ij} )}$ such that $r_{(\lambda_1 , \lambda_2) }(z_1,z_2)=z^{m_i(\lambda_i)m_j(\lambda_j)}=r_{(\lambda_1^{ij} , \lambda_2^{ij}) }(w_1,w_2)$. Set $\omega_m$ be the $m$-th root of 1. Then $(z_1,z_2),\ (w_1,w_2)$ are restrict to the following forms:
\begin{equation*}
z_{p_1,p_2}=(z_1,z_2)=\bigl(z^{m_j(\lambda_2)} \omega_{|m_i(\lambda_1)|}^{p_1}, z^{n_i(\lambda_1)} \omega_{|m_i(\lambda_1)m_j(\lambda_2)|}^{n_i(\lambda_1)p_1 + m_i(\lambda_1)p_2} \bigr)
\end{equation*}
for $0 \le p_1 \le |m_i(\lambda_1)|-1, 0 \le p_2 \le |m_j(\lambda_2)|-1$ and
\begin{equation*}
z'_{q_1,q_2}=(w_1,w_2)=\bigl(z^{m_i(\lambda_2^{ij})} \omega_{|m_j(\lambda_1^{ij})|}^{q_1} , z^{n_j(\lambda_1^{ij})} \omega_{|m_j(\lambda_1^{ij}) m_i(\lambda_2^{ij})|}^{n_j(\lambda_1^{ij})q_1+m_j(\lambda_1^{ij})q_2} \bigr), \  
\end{equation*}
for $0 \le q_1 \le |m_j(\lambda_1^{ij})|-1 , 0 \le q_2 \le |m_i(\lambda_1^{ij})|-1$.
 Each $(p_1,p_2) \in \{ 0 ,\cdots , |m_i(\lambda_1)|-1\} \times \{ 0 ,\cdots , |m_j(\lambda_2)|-1\}$, there exists the unique element $(p_1^{ij},p_2^{ij}) \in \{ 0 ,\cdots , |m_j(\lambda^{ij})|-1\} \times \{ 0 ,\cdots , |m_i(\lambda_2^{ij})|-1\}$ with the relation $n_j(\lambda_2)(n_i(\lambda_1)p_1 + m_i(\lambda_i)p_2)=n_i(\lambda_2^{ij})(n_j(\lambda_1^{ij})p_1^{ij}+m_j(\lambda_1^{ij})p_2^{ij})$ modulo $|m_i(\lambda_1)m_j(\lambda_2)|$ by the assumption (i),(ii).
 Hence if we define $S_{i,j}(z_{p_1,p_2})=z'_{p_1^{ij},p_2^{ij}}$, then this map is homeomorphism and preserving range and source maps.
 Moreover, we can check that $S_{i,j}$ satisfies the hexagonal condition by a tedious calculation.
 Hence we can construct a topological $k$-graph $\Lambda_{(\Gamma ,m ,n)}$ from $\Gamma^{e_i} \times_{n_i,m_i} \mathbb{T}$ and $\widetilde{T}_{i,j}=T_{i,j} \times S_{i,j}$. This topological $k$-graph is row-finite with no source by Proposition \ref{prop:test}.
\end{proof}

\begin{ex}[Cuntz's $ax+b$-semigroup $C^*$-algebra over $\mathbb{N}$]
\label{ex:ax+b}
 Assume $k=\infty$.
 Let us consider the discrete set $\Lambda^0=\{ v_0 \}$ and $\Lambda^{e_i}=\{ \lambda_i \}$. Define $T_{e_i,e_j}: \Gamma^{e_i} \times_c \Gamma^{e_j} \longrightarrow \Gamma^{e_j} \times_c \Gamma^{e_i}$ by $T_{e_i,e_j}(\lambda_i,\lambda_j)=(\lambda_j,\lambda_i)$.
Let $\mathcal{P}=\{ p_1,p_2,\cdots \}$ be the set of prime numbers with $p_1 \le p_2 \le \cdots$.
 Define $m_i(\lambda_i)=1$, $n_i(\lambda_i)=p_i$. Then $\Gamma^{e_i}, m_i, n_i$ satisfy the assuption of Proposition \ref{prop:newgraph}, hence we can construct a topological $k$-graph $\Lambda_{\Gamma , m , n}$. For this graph, the $C^*$-algebra $\mathcal{O}(\Lambda_{\Gamma , m , n})$ is isomorphic to the $ax+b$-semigroup $C^*$-algebra $\mathcal{Q}_\mathbb{N}$ defined in \cite{Cuntz} which is the universal $C^*$-algebra generated by isometries $s_n$ $(n \in \mathbb{N}^*)$ and a unitary $u$ satisfying the relations
\begin{equation*}
s_ns_m=s_{nm}, \quad s_nu=u^ns_n , \quad \sum_{i=0}^{n-1}u^i  s_ns_n^* u^{-i}=1
\end{equation*}
for $n,m \in \mathbb{N}^*$.
 Let $a_0$ be a generator of $A=C(\mathbb{T})$ and $\xi_i=1/\sqrt{p_i}$ be a constant function of $X_{e_i}=C(\mathbb{T})$, then $u$ and $s_{p_i}$ are corresponding to $t^0(a_0)$ and $t^{e_i}( \xi_i)$ respectively.
 
 Cuntz showed that the $C^*$-algebra $\mathcal{Q}_\mathbb{N}$ is simple and purely infinite in \cite{Cuntz}. This property is also followed by Theorem \ref{thm:purelyinf} in our state.
\end{ex}

\begin{ex}
 Let us consider the arbitrary $k=1,2,\cdots ,\infty$. Set $\Gamma$ and $T$ as in Example \ref{ex:ax+b}.
 Let $\{ (p_i,q_i) \}_{i=1}^k$ be a set of pairs such that $p_i \in \mathbb{Z} \setminus \{ 0 \}$ and $q_i \in \mathbb{Z}_+$ with $\mathrm{gcd}(|p_i|,q_j)=1$ $(1 \le i,j \le k)$. Set $m_i(\lambda_i)=p_i$ and $n_i(\lambda_i)=q_i$. By Proposition \ref{prop:newgraph}, we can define a topological $k$-graph $\Lambda_{\Gamma , m , n}$.
 If $(|p_i|,q_i) \neq (1,1)$ for all $i$, then $\Lambda$ satisfies condition (A).
 Under this condition, if there exists $1 \le i \le k$ such that $p_i \notin q_i\mathbb{Z}$, then $\mathcal{O}(\Lambda_{\Gamma , m , n})$ is simple purely infinite by Theorem \ref{thm:purelyinf}.
 In \cite{Cuntz}, Cuntz showed $\mathcal{Q}_\mathbb{N}$ is generated by Bost-Connes algebra (\cite{BostConnes}) adding one unitary. So we have an interesting the corresponding part of Bost-Connes algebra in this example (See Remark \ref{rem:BostConnes}). 
\end{ex}

Next, we give two examples of topological $2$-graphs whose associated (discrete) $2$-graphs are not only one loop.

\begin{ex}
 Let us consider the following (discrete) 2-graph $\Gamma$;
\[
\Gamma=
\xymatrix {
{\bullet} \ar@(dl,ul)^{\lambda_1}[] \ar@/^/@{-->}[rr]^{\mu_2}  &&
{\bullet} \ar@/^/@{-->}[ll]^{\mu_1} \ar@(dr,ur)_{\lambda_2}[] 
}
.\]
$\Gamma$ is characterized by $\Gamma^{e_1}=\{ \lambda_1 , \lambda_2 \}$, $\Gamma^{e_2}=\{ \mu_1, \mu_2 \}$, and $T_{e_1,e_2}(\lambda_1, \mu_1)=(\mu_1, \lambda_2)$, $T_{e_1,e_2}(\lambda_2, \mu_2)=(\mu_2, \lambda_1)$. 
 By (ii) of Proposition \ref{prop:newgraph}, $p_0=m_1(\lambda_1)=m_1(\lambda_2)$ and $q_0=n_1(\lambda_1)=n_2(\lambda_2)$ have to hold. Set $p_i=m_2(\mu_i)$, $q_i=n_2(\mu_i)$ $(i=1,2)$. 
 Then we have to impose $\mathrm{gcd}(|p_0p_i|,q_0q_i)=1$ for $i=1,2$.
 If $(|p_0|,q_0) \neq (1,1)$ and $(|p_1p_2|,q_1q_2) \notin \mathbb{Z}_+(1,1)$, then $\Lambda_{\Gamma , m ,n}$ is Condition (A). Moreover $p_0 \notin q_0 \mathbb{Z}$ or $p_1p_2 \notin q_1q_2 \mathbb{Z}$, then $\mathcal{O}(\Lambda_{\Gamma ,m,n})$ is simple and purely infinite because of the existence of the loop in $\Gamma$. 
\end{ex}

\begin{ex}
 Let us consider the 2-graph $\Omega_2$.
 In a suitable condition, $\mathcal{O}(\Lambda_{\Omega_2 , m , n})$ is a simple $A\mathbb{T}$-algebra (also see the dichotomy after Proposition 3.14 of \cite{Katsura4}).
\end{ex}

\begin{rem}
\label{rem:BostConnes}
If we remove the "no source" assumption, we can treat the Bost-Connes algebra $\mathcal{C}_\mathbb{Q}$ in the framewark of topological $k$-graph.
 Let $\mathcal{P}=\{ p_1,p_2,\cdots \}$ be the set of prime numbers with $p_1 \le p_2 \le \cdots$.
 Set $\mathbb{Z}_p$ ($p \in \mathcal{P}$) be the $p$-adic ring and $\mathcal{Z}=\prod_{p \in \mathcal{P}} \mathbb{Z}_p$. Laca (Proposition 32 of \cite{Laca}) showed $\mathcal{C}_\mathbb{Q}$ is isomorphic to an endomorphism crossed product $C(\mathcal{Z}) \rtimes_\alpha \mathbb{N}^\infty$ where $\alpha_{e_i} : C(\mathcal{Z}) \longrightarrow C(\mathcal{Z})$ is the endomorphism defined by
\[ 
\alpha_{e_i}(f)(x) =
\left\{
  \begin{array}{cl}
  f(x/p_i) & \mbox{if $x \in p_i\mathcal{Z}$} \\
  0 & \mbox{otherwise.}
  \end{array}
\right.
\]
On the other hand, $\Lambda^0=\mathcal{Z}$, $\Lambda^{e_i}=p_i\mathcal{Z}$ $(1 \le i < \infty )$, and define $r_{e_i},\ s_{e_i} : \Lambda^{e_i} \longrightarrow \Lambda^0$ by $r_{e_i}(x)=x$ and $s_{e_i}(x)=x/p_i$, then we can "construct" a topological $\infty$-graph $\Lambda$ such that $\mathcal{O}(\Lambda ) \cong \mathcal{C}_\mathbb{Q}$ (see also Example 2.5 (4) of \cite{Yeend2}). However $\Lambda$ has sources, so far, we cannot treat the Bost-Connes algebra
 $\mathcal{C}_\mathbb{Q}$ in this paper.  
\end{rem}

\end{document}